# Some elementary observations on Narayana polynomials and related topics

Johann Cigler

Fakultät für Mathematik

Universität Wien

johann.cigler@univie.ac.at

**Abstract.**
We give an elementary account of generalized Fibonacci and Lucas polynomials whose moments are Narayana polynomials of type A and type B.

**Introduction**

Consider the Fibonacci polynomials $F_n(x) = \sum_{j=0}^{\lfloor \frac{n}{2} \rfloor} (-1)^j \binom{n-j}{j} x^{n-2j}$ and the corresponding Lucas polynomials $L_n(x) = \sum_{j=0}^{\lfloor \frac{n}{2} \rfloor} (-1)^j \frac{n}{n-j} \binom{n-j}{j} x^{n-2j}$ and let $L$ be the linear functional defined by $L(F_n(x)) = [n=0]$ and $M$ be the linear functional defined by $M(L_n(x)) = [n=0]$. Then the moments $L(x^{2n}) = C_n$ are Catalan numbers and the moments $M(x^{2n}) = M_n = \binom{2n}{n}$ are central binomial coefficients. An analogous situation holds by replacing the Catalan numbers $C_n$ by the Narayana polynomials $C_n(t) = \sum_{k \geq 0} \binom{n-1}{k}\binom{n}{k}\frac{1}{k+1} t^k$ and the central binomial coefficients $M_n$ by the polynomials $M_n(t) = \sum_{j=0}^{n} \binom{n}{j}^2 t^j$, which are sometimes called Narayana polynomials of type B.

In this survey article I give an elementary and self-contained account of the corresponding polynomials and the associated Catalan-Stieltjes matrices. I want to thank Dennis Stanton and Jiang Zeng for helpful remarks and references to the literature.

## 1. 1. Background material on Fibonacci polynomials and Catalan numbers

The basic facts about Fibonacci and Lucas polynomials are very old and well known (cf. e.g. [5]).

The *Fibonacci polynomials* $f_n(x,s) = \sum_{k=0}^{\lfloor \frac{n-1}{2} \rfloor} \binom{n-1-k}{k} x^{n-1-2k} s^k$ satisfy the recursion

$f_n(x,s) = xf_{n-1}(x,s) + sf_{n-2}(x,s)$ with initial values $f_0(x,s) = 0$ and $f_1(x,s) = 1$.

We will consider the *special Fibonacci polynomials* $F_n(x) = f_{n+1}(x,-1)$. If $U_n(x)$ denotes a *Chebyshev polynomial of the second kind* then we can equivalently write $F_n(x) = U_n\left(\dfrac{x}{2}\right)$.

The first terms of the sequence $\left(F_n(x)\right)_{n\geq 0}$ are

1, x, -1 + x², -2 x + x³, 1 - 3 x² + x⁴, 3 x - 4 x³ + x⁵, ...

**Remark**

Let me recall some well-known facts about orthogonal polynomials (cf. [4], [13],[17]). These are polynomials $\left(p_n(x)\right)_{n\geq -1}$ satisfying a recursion of the form
$p_n(x) = (x - s_{n-1})p_{n-1}(x) - t_{n-2}p_{n-2}(x)$ with initial values $p_{-1}(x) = 0$ and $p_0(x) = 1$. The corresponding *Catalan-Stieltjes matrix* $(a(n,k))$ (cf. [13]) consists of the uniquely determined numbers $a(n,k)$ which satisfy $x^n = \sum_{k=0}^{n} a(n,k) p_k(x)$.

It satisfies

$$a(n,k) = a(n-1,k-1) + s_k a(n-1,k) + t_k a(n-1,k+1) \tag{1.1}$$

with $a(0,k) = [k=0]$ and $a(n,-1) = 0$ because

$$\sum_{k=0}^{n} a(n,k) p_k(x) = x \cdot x^{n-1} = \sum_{k=0}^{n} a(n-1,k) x p_k(x) = \sum_{k=0}^{n} a(n-1,k)\left(p_{k+1}(x) + s_k p_k(x) + t_{k-1} p_{k-1}(x)\right)$$
$$= \sum a(n-1,k-1) p_k(x) + \sum s_k a(n-1,k) p_k(x) + \sum t_k a(n-1,k+1) p_k(x).$$

The numbers $s_k$ and $t_k$ uniquely determine both the polynomials $p_n(x)$ and the corresponding Catalan-Stieltjes matrix.

Let $L$ be the linear functional defined by $L(p_n) = [n=0]$. Here we use Iverson's convention $[P] = 1$ if property $P$ is true and $[P] = 0$ else. The polynomials satisfy moreover $L(p_n p_m) = 0$ for $m \neq n$, i.e. they are orthogonal with respect to $L$. But we shall not use this property.

The numbers $L(x^n)$ are called moments of the sequence $(p_n(x))$.

If all $s_k = 0$ then $P_n(x) = p_{2n}(\sqrt{x})$ satisfies

$P_1(x) = x - t_0$ and $P_n(x) = (x - t_{2n-1} - t_{2n})P_{n-1}(x) - t_{2n} t_{2n+1} P_{n-2}(x)$

and $Q_n(x) = \dfrac{p_{2n+1}(\sqrt{x})}{\sqrt{x}}$ satisfies $Q_n(x) = (x - t_{2n} - t_{2n+1})Q_{n-1}(x) - t_{2n+1} t_{2n+2} Q_{n-2}(x).$

This splitting is equivalent with the odd-even trick in [6].

For the Fibonacci polynomials $F_n(x)$ the numbers $a(n,k)$ satisfy

$$a(n,k) = a(n-1, k-1) + a(n-1, k+1) \tag{1.2}$$

with $a(0,k) = [k=0]$.

Thus $a(n,k)$ can be interpreted as the number of elements of the set of $n-$ letter words $w_1 w_2 \cdots w_n$ in the alphabet $\{-1, 1\}$ that add up to $k$, and all whose partial sums are non-negative because for $w_n = 1$ the word $w_1 w_2 \cdots w_{n-1}$ adds up to $k-1$ and for $w_n = -1$ to $k+1$.

These so-called *ballot numbers* are well known and satisfy

$$a(2n+k, k) = \binom{2n+k}{n} - \binom{2n+k}{n-1}. \tag{1.3}$$

or equivalently

$$x^n = \sum_{k=0}^{\lfloor \frac{n}{2} \rfloor} \left( \binom{n}{k} - \binom{n}{k-1} \right) F_{n-2k}(x). \tag{1.4}$$

Let $L$ be the linear functional defined by $L(F_n) = [n=0]$. Here $[P] = 1$ if property $P$ is true and $[P] = 0$ else. Then (1.4) implies

$$L(x^{2n}) = \binom{2n}{n} - \binom{2n}{n-1} = C_n = \binom{2n}{n} \frac{1}{n+1} \tag{1.5}$$

is a *Catalan number* and $L(x^{2n+1}) = 0$.

The first terms of the sequence $(C_n)_{n \geq 0}$ are

1, 1, 2, 5, 14, 42, 132, 429, 1430, 4862, ...

Let us compute the generating functions $f_k(z) = \sum_{n \geq 0} a(n,k) z^n$. Then (1.2) translates into

$$f_k(z) = z \left( f_{k-1}(z) + f_{k+1}(z) \right) \tag{1.6}$$

and

$$f_0(z) = 1 + z f_1(z). \tag{1.7}$$

The uniquely determined solution of these equations is $f_k(z) = z^k f(z)^{k+1}$ if we set $f(z) = f_0(z)$.

This can easily be verified by comparing coefficients.

By (1.7) $f(z)$ satisfies $f(z) = 1 + z^2 f(z)^2$ which implies the well-known result

$$f(z) = \sum_{n \geq 0} C_n z^{2n} = \frac{1 - \sqrt{1-4z^2}}{2z^2}. \tag{1.8}$$

Let us also consider the polynomials

$$P_n(x) = F_{2n}\left(\sqrt{x}\right) = \sum_{k=0}^{n} (-1)^{n-k} \binom{n+k}{2k} x^k \tag{1.9}$$

and

$$Q_n(x) = \frac{F_{2n+1}\left(\sqrt{x}\right)}{\sqrt{x}} = \sum_{k=0}^{n} (-1)^{n-k} \binom{n+k+1}{2k+1} x^k. \tag{1.10}$$

By (1.4) we get

$$x^n = \sum_{k=0}^{n} \left(\binom{2n}{k} - \binom{2n}{k-1}\right) P_{n-k}(x). \tag{1.11}$$

Let $L_0$ denote the linear functional defined by $L_0(P_n) = [n=0]$.

Then we get for the moments

$$L_0(x^n) = C_n. \tag{1.12}$$

Analogously we get

$$x^n = \sum_{k=0}^{n} \left(\binom{2n+1}{k} - \binom{2n+1}{k-1}\right) Q_{n-k}(x). \tag{1.13}$$

Let $L_1$ denote the linear functional defined by $L_1(Q_n) = [n=0]$.

Then we get for the moments

$$L_1(x^n) = \binom{2n+1}{n} - \binom{2n+1}{n-1} = C_{n+1}. \tag{1.14}$$

### 1.2. Narayana polynomials as moments

The Catalan numbers are special cases for $t=1$ of the *Narayana polynomials*

$$C_n(t) = \sum_{k \geq 0} \binom{n-1}{k} \binom{n}{k} \frac{1}{k+1} t^k \tag{1.15}$$

for $n > 0$ and $C_0(t) = 1$. (cf. [14]).

The first terms of $(C_n(t))_{n\geq 0}$ are

1, 1, 1 + t, 1 + 3 t + t², 1 + 6 t + 6 t² + t³, 1 + 10 t + 20 t² + 10 t³ + t⁴, ...

For $t = 2$ they reduce to the *little Schroeder numbers* $(C_n(2))_{n\geq 0} = (1,1,3,11,45,197,\cdots)$, OEIS [12], A001003.

Let $\tau_{2n}(t) = 1$ and $\tau_{2n+1}(t) = t$. Define polynomials $F_n(x,t)$ by the recursion

$$F_n(x,t) = xF_{n-1}(x,t) - \tau_{n-2}(t)F_{n-2}(x,t) \quad (1.16)$$

with initial values $F_0(x,t) = 1$ and $F_1(x,t) = x$.

The first terms of the sequence $(F_n(x,t))_{n\geq 0}$ are

1, x, -1 + x², -x - t x + x³, 1 - 2 x² - t x² + x⁴, x + t x + t² x - 2 x³ - 2 t x³ + x⁵, ...

Their generating function is

$$\sum_{n\geq 0} F_n(x,t)z^n = \frac{1 + xz + tz^2}{1 - (x^2 - 1 - t)z^2 + tz^4}. \quad (1.17)$$

Then we get

**Theorem 1 ([1],[3], [11], [13], [16],[17])**

Let $L$ be the linear functional defined by $L(F_n(x,t)) = [n = 0]$. Then the moments satisfy

$$\begin{aligned} L(x^{2n}) &= C_n(t), \\ L(x^{2n+1}) &= 0. \end{aligned} \quad (1.18)$$

**Remark**

By starting with $C_n(t)$ it is easy to guess (1.16) in the same manner as I have done in [4].

In order to guess explicit formulae for $F_n(x,t)$ it is convenient to consider the polynomials with odd and even degrees separately. To this end we consider the polynomials

$$P_n(x,t) = F_{2n}(\sqrt{x},t) \text{ and } Q_n(x,t) = \frac{F_{2n+1}(\sqrt{x},t)}{\sqrt{x}}.$$

Then (1.32) and (1.21) can be summarized to give the formula

$$F_n(x,t) = \sum_{k=0}^{\lfloor \frac{n}{2} \rfloor} (-1)^k \sum_{j=0}^{k} \binom{\lfloor \frac{n}{2} \rfloor - j}{k - j}\binom{\lfloor \frac{n-1}{2} \rfloor - k + j}{j} t^j x^{n-2k}. \quad (1.19)$$

**1.2.1.** The polynomials $Q_n(x,t)$.

The polynomials $Q_n(x,t)$ satisfy the recurrence

$$Q_n(x,t) = (x-1-t)Q_{n-1}(x,t) - tQ_{n-2}(x,t) \tag{1.20}$$

with initial values $Q_0(x,t) = 1$ and $Q_1(x,t) = x-1-t$.

Thus $Q_n(x,t) = f_{n+1}(x-1-t,-t)$. Binet's formula gives $Q_n(x,t) = \dfrac{\alpha^{n+1} - \beta^{n+1}}{\alpha - \beta}$

with $\alpha = \alpha(x,t) = \dfrac{x-1-t+\sqrt{(x-1-t)^2 - 4t}}{2}$ and $\beta = \beta(x,t) = \dfrac{x-1-t-\sqrt{(x-1-t)^2 - 4t}}{2}$.

A more general class of polynomials has been considered in [1].

By induction we get $Q_n(x,t) = \sum\limits_{k=0}^{n}(-1)^{n-k} q_{n,k}(t) x^k$ with

$$q_{n,k}(t) = \sum_{j=0}^{n-k}\binom{n-j}{k}\binom{k+j}{j}t^j. \tag{1.21}$$

From (1.10) we see that $q_{n,k}(1) = \binom{n+k+1}{2k+1}$.

The first terms of $q_{n,k}(t)$ are

```
1
1 + t                1
1 + t + t²           2 + 2 t              1
1 + t + t² + t³      3 + 4 t + 3 t²      3 + 3 t             1
1 + t + t² + t³ + t⁴ 4 + 6 t + 6 t² + 4 t³  6 + 9 t + 6 t²   4 + 4 t           1
1 + t + t² + t³ + t⁴ + t⁵ 5 + 8 t + 9 t² + 8 t³ + 5 t⁴ 10 + 18 t + 18 t² + 10 t³  10 + 16 t + 10 t²  5 + 5 t  1
```

Note that the polynomials $q_{n,k}(t)$ are palindromic.

Let $B_{n,k}(t)$ be the uniquely determined polynomials such that

$$x^n = \sum_{k=0}^{n} B_{n,k}(t) Q_k(x,t). \tag{1.22}$$

The recursion of $Q_n(x,t)$ implies that

$$B_{n,k}(t) = B_{n-1,k-1}(t) + (1+t)B_{n-1,k}(t) + tB_{n-1,k+1}(t) \tag{1.23}$$

with $B_{0,k}(t) = [k=0]$ and $B_{n,-1}(t) = 0$.

The first terms of the sequence $\bigl(B_{n,0}(t), B_{n,1}(t), \cdots, B_{n,n}(t)\bigr)_{n \geq 0}$ are

```
1
1 + t                        1
1 + 3 t + t²                 2 + 2 t              1
1 + 6 t + 6 t² + t³          3 + 8 t + 3 t²       3 + 3 t           1
1 + 10 t + 20 t² + 10 t³ + t⁴  4 + 20 t + 20 t² + 4 t³   6 + 15 t + 6 t²   4 + 4 t   1
```

By induction we can verify that

$$B_{n,k}(t) = \sum_{j=0}^{n}\binom{n+1}{k+1+j}\binom{n+1}{j}\frac{k+1}{n+1}t^j = \sum_j \left(\binom{n}{j}\binom{n+1}{k+j+1} - \binom{n+1}{j}\binom{n}{k+j+1}\right)t^j. \quad (1.24)$$

For $k=0$ we get

$$B_{n,0}(t) = C_{n+1}(t). \quad (1.25)$$

From (1.13) we see that $B_{n,k}(1) = \binom{2n+1}{n-k} - \binom{2n+1}{n-k-1} = \frac{2k+2}{n+k+2}\binom{2n+1}{n-k}.$

This gives the Catalan triangle OEIS[12], A039598

$$\begin{pmatrix} 1 & 0 & 0 & 0 & 0 \\ 2 & 1 & 0 & 0 & 0 \\ 5 & 4 & 1 & 0 & 0 \\ 14 & 14 & 6 & 1 & 0 \\ 42 & 48 & 27 & 8 & 1 \end{pmatrix}$$

For the little Schroeder numbers the corresponding triangle is OEIS [12], A110440,

$$\begin{pmatrix} 1 & 0 & 0 & 0 & 0 \\ 3 & 1 & 0 & 0 & 0 \\ 11 & 6 & 1 & 0 & 0 \\ 45 & 31 & 9 & 1 & 0 \\ 197 & 156 & 60 & 12 & 1 \end{pmatrix}$$

There is a nice interpretation in terms of weighted NSEW-paths. A *NSEW-path* is a path consisting of North, South, East and West steps of length 1. (Cf. [9] and [10]). We consider only NSEW- paths which start at $(0,0)$ and end on height $k \geq 0$ and never cross the $x$ – axis.

$B_{n,k}(t)$ is the weight of all those NSEW-paths with $n$ steps which end on height $k$, if the weight is defined by $w(N) = w(E) = 1$ and $w(S) = w(W) = t$. This follows immediately from (1.23) because there are 4 possibilities to reach a point of height $k$. For $k=0$ this reduces to

$B_{n,0}(t) = (1+t)B_{n-1,0}(t) + tB_{n-1,1}(t).$

For example for $n=2$ and $k=0$ we get $w(EE) = 1, w(NS + EW + WE) = 3t, w(WW) = t^2$.

For $k=1$ we get $w(NE)+w(EN)=2$ and $w(NW)+w(WN)=2t$.

Let $y \geq 0$ and let $w_n(x,y)$ be the number of NSEW-paths from $(0,0)$ to $(x,y)$ which do not cross the $x-$ axis. It has been shown in [9] that

$$w_n(-n+k+2j,k) = \binom{n}{j}\binom{n}{k+j} - \binom{n}{j-1}\binom{n}{k+j+1} = \binom{n+1}{k+1+j}\binom{n+1}{j}\frac{k+1}{n+1}.$$

A purely combinatorial proof has been given in [10] and can be considered as another proof of (1.24).

All these polynomials are palindromic and *gamma-nonnegative*, i.e. they have a representation of the form $\sum \gamma_{n,j} t^j (1+t)^{n-2j}$ where $\gamma_{n,j}$ are non-negative integers. (Cf. [14] for this notion).

More precisely we have

$$B_{n,k}(t) = \sum_{i=0}^{\lfloor \frac{n-k}{2} \rfloor} \binom{k+2i}{i}\frac{k+1}{i+k+1}\binom{n}{2i+k} t^i (1+t)^{n-k-2i}, \tag{1.26}$$

which for $k=0$ reduces to

$$C_{n+1}(t) = \sum_{i=0}^{\lfloor \frac{n}{2} \rfloor} C_i \binom{n}{2i} t^i (1+t)^{n-2i}. \tag{1.27}$$

In order to prove this we modify a method developed in [15]. Let $f(N)=1, f(S)=-1, f(E)=f(W)=0$.

To each non-negative NSEW- path $u_1 \cdots u_n$ with $u_i \in \{N,S,E,W\}$ whose endpoint is on height $k$ we associate the $n-$ letter word $f(u_1)f(u_2)\cdots f(u_n)$ in the alphabet $\{-1,1,0\}$ that adds up to $k$, and all whose partial sums are non-negative.

For each such sequence there are $i$ terms $f(u_j)=-1$ and $i+k$ terms $f(u_j)=1$ for some $i$.

On the other hand we can choose $2i+k$ places where $u_j = N$ or $u_j = S$, i.e. $f(u_j) = \pm 1$ in $\binom{n}{2i+k}$ ways. By (1.3) we can order the signs in such a way that the corresponding path is non-negative in $\binom{k+2i}{i} - \binom{k+2i}{i-1} = \binom{k+2i}{i}\frac{k+1}{i+k+1}$ ways. In the remaining $n-2i-k$ places we can arbitrarily put $W$ or $E$. The weight of all such paths is therefore
$\binom{n}{2i+k}\binom{k+2i}{i}\frac{k+1}{i+k+1} t^i (1+t)^{n-k-2i}$.

If we define the linear functional $L_1$ by $L_1(Q_n(x,t)) = [n=0]$ we get from (1.27) that

$$L_1(x^n) = C_{n+1}(t). \tag{1.28}$$

Let us compute the generating functions $f_k(z,t) = \sum_{n\geq 0} B_{n,k}(t) z^n$. As above we see that they satisfy

$$f_k(z,t) = z\big(f_{k-1}(z,t) + (1+t)f_k(z,t) + tf_{k+1}(z,t)\big) \text{ with } f_0(z,t) = 1 + (1+t)zf_0(z,t) + tzf_1(z,t).$$

The unique solution is

$$f_k(z,t) = z^k f(z,t)^{k+1} \text{ where } f(z,t) \text{ satisfies } 1 - (1-(1+t)z)f(z,t) + tz^2 f(z,t)^2 = 0.$$

This implies

$$f(z,t) = \sum_{n\geq 0} C_{n+1}(t) z^n = \frac{1-(1+t)z - \sqrt{(1-(1+t)z)^2 - 4tz^2}}{2tz^2}. \tag{1.29}$$

Since $1 - (1-(1+t)z)f(z,t) + tz^2 f(z,t)^2 = 0$ we get

$$\sum_k B_{n,k}(t) \frac{t^{k+1}-1}{t-1} z^n = \frac{1}{t-1}\left( \sum_k z^k f(z,t)^{k+1} t^{k+1} - \sum_k z^k f(z,t)^{k+1} \right) = \frac{f(z,t)}{t-1}\left( \frac{t}{1-tzf(z,t)} - \frac{1}{1-zf(z,t)} \right)$$

$$= \frac{f(z,t)}{t-1} \frac{(t-1)}{(1-zf(z,t))(1-tzf(z,t))} = \frac{f(z,t)}{1-(1+t)zf(z,t) + tz^2 f(z,t)^2} = \frac{f(z,t)}{f(z,t) - 2(1+t)zf(z,t)} = \frac{1}{1-2(1+t)z}.$$

This implies

$$\sum_{k=0}^n B_{n,k}(t)\big(1+t+\cdots+t^k\big) = (2t+2)^n. \tag{1.30}$$

A combinatorial proof of (1.30) has been given in [2], proof of identity 1, in a somewhat different context which we will translate into our terminology.

The right-hand side of (1.30) is the weight of all NSWE-paths of length $n$.

Let $\mathbf{B}_{n,k}$ be the set of all non-negative NSWE-paths of length $n$ which end on height $k$.

For $p \in \mathbf{B}_{n,k}$ we define $k+1$ different paths $\varphi_i(p)$, $0 \leq i \leq k$, of length $n$ such that $w(\varphi_i(p)) = t^i w(p)$.

To this end define the last ascent to height $i$ of $p$ to be the last step $N$ from height $i-1$ to $i$. Let $\varphi_i(p)$ denote the path obtained by changing each of the last ascents to heights $1, 2, \cdots, i$ to downsteps $S$. For $i = 0$ let $\varphi_0(p) = p$. Then all $\varphi_i(p)$ are different and for $i > 0$ not non-negative. The height of $\varphi_i(p)$ is $k - 2i$ and the weight is $w(\varphi_i(p)) = t^i w(p)$.

Let on the other hand $q$ be a path with height $j$, which crosses the $x-$ axis. Then it has a set of premier descents below the $x-$ axis, i.e. the first (from left to right) down steps $S$ from height $m$ to $m-1$ for $m=0,-1,\cdots$. Suppose $q$ has $i$ premier descents below the $x-$ axis. Then changing each of these $S$ to upsteps $N$ gives a new path $p$ which is non-negative and ends on height $j+2i$. It is clear that $\varphi_i(p)=q$ and $w(\varphi_i(p))=t^i w(p)$.

For example

$\mathbf{B}(2,0)=\{EE,EW,WE,NS,WW\}$,

$\mathbf{B}(2,1)=\{NE,NW,EN,WN\}$, $\varphi_1(\mathbf{B}(2,1))=\{SE,SW,ES,WS\}$,

$\mathbf{B}(2,2)=\{NN\}$, $\varphi_1(\mathbf{B}(2,2))=\{SN\}$, $\varphi_2(\mathbf{B}(2,2))=\{SS\}$.

### 1.2.2. The polynomials $P_n(x,t)$.

The polynomials $P_n(x,t)$ satisfy the recurrence

$$P_n(x,t)=(x-\sigma_{n-1}(t))P_{n-1}(x,t)-tP_{n-2}(x,t)$$

with initial values $P_0(x,t)=1$ and $P_1(x,t)=x-1$,

where $\sigma_0(t)=1$ and $\sigma_n(t)=1+t$ for $n>0$.

We have for $n>0$

$$P_n(x,t)=Q_n(x,t)+tQ_{n-1}(x,t). \tag{1.31}$$

For (1.31) holds for $n=1$ and $n=2$ and for $n\geq 3$ both sides satisfy the same recursion.

Let us set $P_n(x,t)=\sum_{k=0}^{n}(-1)^{n-k}p_{n,k}(t)x^k$.

Then we get

$$p_{n,k}(t)=\sum_{j=0}^{n-k}\binom{n-j}{k}\binom{k-1+j}{j}t^j. \tag{1.32}$$

The first terms of the sequence
$\left(p_{n,0}(t)=q_{n,0}(t)-q_{n-1,0}(t),\, p_{n,1}(t)=q_{n,1}(t)-q_{n-1,1}(t),\cdots,p_{n,n}(t)=q_{n,n}(t)-q_{n-1,n}(t)\right)_{n\geq 0}$ are

```
1
1    1
1    2 + t                      1
1    3 + 2 t + t²               3 + 2 t              1
1    4 + 3 t + 2 t² + t³        6 + 6 t + 3 t²      4 + 3 t          1
1    5 + 4 t + 3 t² + 2 t³ + t⁴ 10 + 12 t + 9 t² + 4 t³  10 + 12 t + 6 t²  5 + 4 t   1
```

Let $A_{n,k}(t)$ be the uniquely determined polynomials satisfying

$$x^n = \sum_{k=0}^{n} A_{n,k}(t) P_k(x,t). \tag{1.33}$$

Then

$$A_{n,k}(t) = A_{n-1,k-1}(t) + \sigma_k(t) A_{n-1,k}(t) + t A_{n-1,k+1}(t) \tag{1.34}$$

with $A_{0,k}(t) = [k=0]$ and $A_{n,-1}(t) = 0$.

This means that $A_{n,k}(t)$ can be interpreted as the weight of all NSEW - paths of length $n$ which end on height $k$ and which have no W-step on height $0$.

For example let $n=3$. For $k=0$ we have $w(EEE)=1$, $w(NSE+ENS+NES)=3t$ and $w(NWS)=t^2$. For $k=2$ we have $w(NNE+ENN+NEN)=3$ and $w(NNW+NWN)=2t$.

The first terms of the sequence $\left(A_{n,0}(t), A_{n,1}(t), \cdots, A_{n,n}(t)\right)_{n \geq 0}$ are

```
1
1                       1
1 + t                   2 + t                   1
1 + 3 t + t²            3 + 5 t + t²            3 + 2 t             1
1 + 6 t + 6 t² + t³     4 + 14 t + 9 t² + t³    6 + 11 t + 3 t²     4 + 3 t    1
```

From (1.31) we get $A_{n,k} + t A_{n,k+1} = B_{n,k}$.

In general we get for $n > 0$

$$A_{n,k}(t) = \sum_{j=0}^{n-k} \binom{n-1}{j}\binom{n}{k+j} \frac{kn+n-j}{(n-j)(k+1+j)} t^j$$
$$= \sum_{j=0}^{n-k} \left(\binom{n-1}{j}\binom{n+1}{k+j+1} - \binom{n}{j}\binom{n}{k+j+1}\right) t^j. \tag{1.35}$$

For $k=0$ this reduces to

$$A_{n,0}(t) = C_n(t). \tag{1.36}$$

For $t=1$ we get the triangle OEIS [12], A039599,

$$\begin{pmatrix} 1 & 0 & 0 & 0 & 0 \\ 1 & 1 & 0 & 0 & 0 \\ 2 & 3 & 1 & 0 & 0 \\ 5 & 9 & 5 & 1 & 0 \\ 14 & 28 & 20 & 7 & 1 \end{pmatrix}$$

For $t=2$ we get OEIS [12], 172094,

$$\begin{pmatrix} 1 & 0 & 0 & 0 & 0 \\ 1 & 1 & 0 & 0 & 0 \\ 3 & 4 & 1 & 0 & 0 \\ 11 & 17 & 7 & 1 & 0 \\ 45 & 76 & 40 & 10 & 1 \end{pmatrix}$$

From (1.33) we get

$$\sum_{k=0}^{n} A_{n,k}(t) F_{2k}(x,t) = x^{2n}. \tag{1.37}$$

Applying the linear functional $L$ gives

$$L(x^{2n}) = A_{n,0}(t) = C_n(t). \tag{1.38}$$

By (1.22) we get $x^{2n+1} = \sum_{k=0}^{n} B_{n,k}(t) F_{2k+1}(x,t)$ which implies $L(x^{2n+1}) = 0$ and thus proves Theorem 1.

If we define the linear functional $L_0$ by $L_0(P_n(x,t)) = [n=0]$ then we get

$$L_0(x^n) = C_n(t). \tag{1.39}$$

Let us also compute the generating functions $f_k(z,t) = \sum_{n \geq 0} A_{n,k}(t) z^n$. They satisfy

$$\begin{aligned} f_k(z,t) &= z\big(f_{k-1}(z,t) + (1+t)f_k(z,t) + tf_{k+1}(z,t)\big), \\ f_0(z,t) &= 1 + z\big(f_0(z,t) + tf_1(z,t)\big). \end{aligned} \tag{1.40}$$

Let $f(z,t)$ satisfy $f(z,t) = 1 + (1+t)zf(z,t) + tz^2 f(z,t)^2$. Then $f_k(z,t) = z^k f_0(z,t) f(z,t)^k$ satisfies the first equation in (1.40). From the second equation and (1.29) we get the well-known formula (cf. e.g. [14])

$$f_0(z) = C(t,z) = \sum_{n \geq 0} C_n(t) z^n = \frac{1 + z(t-1) - \sqrt{1 - 2z(t+1) + z^2(t-1)^2}}{2tz}. \tag{1.41}$$

**Remarks**

In terms of $C(t,z)$ we get

$$\sum_{n\geq 0} A_{n,k}(t)z^n = C(t,z)\big(C(t,z)-1\big)^k,$$

$$\sum_{n\geq 0} B_{n,k}(t)z^n = \frac{\big(C(t,z)-1\big)^{k+1}}{z}.$$

(1.42)

For $t=1$ it is well known that $\big(F_n(1,1)\big) = (1,1,0,-1,-1,0,1,1,0,-1,-1,0,\cdots)$ is periodic with period 6 because $\alpha(1,1) = \dfrac{-1+\sqrt{-3}}{2}$ and $\beta(1,1) = \dfrac{-1-\sqrt{-3}}{2}$ satisfy $\alpha(1,1)^3 = \beta(1,1)^3 = 1$.

For $t=2$ and $t=3$ an analogous situation obtains: $\alpha(1,2) = -1+i$ and $\beta(1,2) = -1-i$ satisfy $\alpha(1,2)^8 = \beta(1,2)^8 = 2^4$ and $\alpha(1,3) = \dfrac{-3+\sqrt{-3}}{2}$ and $\beta(1,3) = \dfrac{-3-\sqrt{-3}}{2}$ satisfy $\alpha(1,3)^{12} = \beta(1,3)^{12} = 3^6$. This implies that the sequence $\left(\dfrac{F_n(1,2)}{4^{\lfloor \frac{n}{8} \rfloor}}\right)_{n\geq 0}$ is periodic with period 16 and the sequence $\left(\dfrac{F_n(1,3)}{27^{\lfloor \frac{n}{12} \rfloor}}\right)_{n\geq 0}$ is periodic with period 24.

We get $\left(\dfrac{F_n(1,2)}{4^{\lfloor \frac{n}{8} \rfloor}}\right)_{n\geq 0} = (1,1,0,-2,-2,2,4,0,-1,-1,0,2,2,-2,-4,0,\cdots)$

and

$\left(\dfrac{F_n(1,3)}{27^{\lfloor \frac{n}{12} \rfloor}}\right)_{n\geq 0} = (1,1,0,-3,-3,6,9,-9,-18,9,27,0,-1,-1,0,3,3,-6,-9,9,18,-9,-27,0,\cdots).$

## 2.1. Background material on Lucas polynomials and central binomial coefficients

The *Lucas polynomials* $l_n(x,s) = \sum_{k=0}^{\lfloor \frac{n}{2} \rfloor} \dfrac{n}{n-k}\binom{n-k}{k} s^k x^{n-2k}$ satisfy the recurrence relation

$l_n(x,s) = x l_{n-1}(x) + s l_{n-2}(x)$ with initial values $l_0(x,s) = 2$ and $l_1(x,s) = x$.

Let us consider the *special Lucas polynomials* $L_n(x)$ defined by $L_n(x) = l_n(x, -1)$ for $n > 0$ and $L_0(x) = 1$.

Then $L_n(x)$ satisfies the recursion

$$L_n(x) = xL_{n-1}(x) - \tau_{n-2}L_{n-2}(x) \tag{2.1}$$

with $\tau_0 = 2$ and $\tau_n = 1$ for $n > 0$.

The first terms of $(L_n(x))_{n \geq 0}$ are

1, x, -2 + x², -3 x + x³, 2 - 4 x² + x⁴, 5 x - 5 x³ + x⁵, ...

Note that $L_n(x) = 2T_n\left(\dfrac{x}{2}\right)$ for $n > 0$ if $T_n(x)$ is a *Chebyshev polynomial of the first kind*.

Let $(a(n,k))$ be the corresponding Catalan-Stieltjes matrix.

Then we get

$a(n,k) = a(n-1, k-1) + a(n-1, k+1)$ for $k > 0$ and $a(n,0) = 2a(n-1,1)$.

Thus $a(n,k)$ is the weight of all non-negative NSEW-paths of length $n$ whose endpoints are on height $k$ where all weights $w(E) = w(N) = w(W) = w(S) = 1$ except that $w(S) = 2$ if the endpoint of $S$ is on the $x-$axis.

The first terms are OEIS [12], A 108044,

$$\begin{pmatrix} 1 & 0 & 0 & 0 & 0 & 0 & 0 \\ 0 & 1 & 0 & 0 & 0 & 0 & 0 \\ 2 & 0 & 1 & 0 & 0 & 0 & 0 \\ 0 & 3 & 0 & 1 & 0 & 0 & 0 \\ 6 & 0 & 4 & 0 & 1 & 0 & 0 \\ 0 & 10 & 0 & 5 & 0 & 1 & 0 \\ 20 & 0 & 15 & 0 & 6 & 0 & 1 \end{pmatrix}$$

This gives $a(2n, 2k) = \dbinom{2n}{n-k}$ and $a(2n+1, 2k+1) = \dbinom{2n+1}{n-k}$ and all other terms vanish.

With other words we get the identities

$$\sum_{k=0}^{n}\binom{2n}{n-k}L_{2k}(x) = x^{2n},$$

$$\sum_{k=0}^{n}\binom{2n+1}{n-k}L_{2k+1}(x) = x^{2n+1}.$$

(2.2)

Let $M$ be the linear functional defined by $M(L_n) = [n=0]$. Then

$$M(x^{2n}) = \binom{2n}{n} \qquad (2.3)$$

is a central binomial coefficient and $M(x^{2n+1}) = 0$.

Let now $f_k(z) = \sum_{n \geq 0} a(n,k) z^n$. Then we have $f_k(z) = f_{k-1}(z) + f_{k+1}(z)$ for $k > 0$ and $f_0(z) = 1 + 2z f_1(z)$. Then we get $f_k(z) = z^k f_0(z) f(z)^k$ with $f(z) = \sum_{n \geq 0} C_n z^n = \dfrac{1 - \sqrt{1-4z}}{2z}$ by (1.8). This gives $f_0(z) = 1 + 2z f_0(z) f(z)$ or

$$f_0(z) = M(z) = \sum_{n \geq 0} \binom{2n}{n} z^n = \frac{1}{\sqrt{1-4z}}. \qquad (2.4)$$

Let us also consider the polynomials

$$R_n(x) = L_{2n}(\sqrt{x}) = \sum_{k=0}^{n}(-1)^{n-k}\frac{2n}{n+k}\binom{n+k}{2k}x^k \qquad (2.5)$$

and

$$S_n(x) = \frac{L_{2n+1}(\sqrt{x})}{\sqrt{x}} = \sum_{k=0}^{n}(-1)^{n-k}\frac{2n+1}{n+k+1}\binom{n+k+1}{2k+1}x^k. \qquad (2.6)$$

Let $M_0$ be the linear functional defined by $M_0(R_n) = [n=0]$. Then (2.2) gives

$$M_0(x^n) = \binom{2n}{n} = M_n. \qquad (2.7)$$

If $M_1$ is the linear functional defined by $M_1(S_n) = [n=0]$ then we get

$$M_1(x^n) = \binom{2n+1}{n} = \frac{1}{2}\binom{2n+2}{n+1} = \frac{M_{n+1}}{2}. \qquad (2.8)$$

## 2.2. The Narayana polynomials of type B as moments

The *central binomial coefficients* are the special case for $t=1$ of the *Narayana polynomials*
$$M_n(t) = \sum_{k=0}^{n} \binom{n}{k}^2 t^k \text{ of type B.}$$

For $t=2$ we get the *central Delannoy numbers* $\left(M_n(2)\right)_{n\geq 0} = (1,3,13,63,321,1683,\cdots)$. Here
$$M_n(2) = d_n = \sum_{k=0}^{n} \binom{2k}{k}\binom{n+k}{2k} = \sum_{k=0}^{n} \binom{n}{k}\binom{n+k}{k}.$$

Let
$$\tau_0(t) = 1+t,$$
$$\tau_{2n}(t) = \frac{1+t^{n+1}}{1+t^n} \text{ for } n>0, \tag{2.9}$$
$$\tau_{2n+1}(t) = \frac{t(1+t^n)}{1+t^{n+1}}.$$

Thus the sequence $\tau_n(t)$ satisfies $\tau_{2n}(t) = 1+t-\tau_{2n-1}(t)$ and $\tau_{2n+1}(t) = \dfrac{t}{\tau_{2n}(t)}$ with initial values $\tau_0(t) = 1+t$ and $\tau_1(t) = \dfrac{2t}{1+t}$.

Define polynomials $L_n(x,t)$ by the recurrence
$$L_n(x,t) = xL_{n-1}(x,t) - \tau_{n-2}(t)L_{n-2}(x,t) \tag{2.10}$$

with initial values $L_0(x,t) = 1$ and $L_1(x,t) = x$.

The first terms of the sequence $\left(L_n(x,t)\right)_{n\geq 0}$ are

$$1,\ x,\ -1-t+x^2,\ -\frac{x(1+4t+t^2-x^2-tx^2)}{1+t},\ 1+t^2-2x^2-2tx^2+x^4,\ \ldots$$

It is clear that $L_n(x,1) = L_n(x)$.

Let now
$$R_n(x,t) = L_{2n}\left(\sqrt{x},t\right). \tag{2.11}$$

These polynomials satisfy
$$R_n(x,t) = (x-1-t)R_{n-1}(x,t) - T_{n-2}(t)R_{n-2}(x,t) \tag{2.12}$$

with

$$T_n(t) = t \text{ for } n > 0,$$
$$T_0(t) = 2t.$$
(2.13)

Then we get
$$R_n(x,t) = Q_n(x,t) - tQ_{n-2}(x,t) \tag{2.14}$$

for $n \geq 2$ and $R_0(x,t) = 1$ and $R_1(x,t) = x - 1 - t$.

For $n > 0$ we get

$$R_n(x,t) = (-1)^n (1+t^n) + \sum_{\ell=1}^{n} (-1)^{n-\ell} \binom{n}{\ell} x^\ell \sum_{j=0}^{n-\ell} \binom{n-\ell}{j} \frac{\binom{\ell+j-1}{j}}{\binom{n-1}{j}} t^j. \tag{2.15}$$

We also have $R_n(x,t) = \alpha^n + \beta^n$ for $n > 0$. This means that $R_n(x,t)$ are the Lucas polynomials corresponding to $Q_n(x,t)$.

If we set $R_0(x,t) = 2$ then the sequence $(R_n(1,1))_{n \geq 0} = (2,-1,-1,\cdots)$ is periodic with period 3, the sequence $\left(\frac{R_n(1,2)}{(2^4)^{\lfloor n/8 \rfloor}}\right)_{n \geq 0} = (2,-2,0,4,-8,8,0,-16,\cdots)$ is periodic with period 8, and the sequence $\left(\frac{R_n(1,3)}{(3^6)^{\lfloor n/12 \rfloor}}\right)_{n \geq 0} = (2,-3,3,0,-9,27,-54,81,-81,0,243,-729,\cdots)$ is periodic with period 12.

Let $D_{n,k}(t)$ be the uniquely determined polynomials such that

$$x^n = \sum_{k=0}^{n} D_{n,k}(t) R_k(x,t). \tag{2.16}$$

They satisfy
$$D_{n,k}(t) = D_{n-1,k-1}(t) + (1+t)D_{n-1,k}(t) + T_k(t)D_{n-1,k+1}(t) \tag{2.17}$$

with $D_{0,k}(t) = [k=0]$ and $D_{n,-1}(t) = 0$.

This implies that
$$D_{n,k}(t) = \left[x^{n-k}\right]\left(1 + (1+t)x + tx^2\right)^n. \tag{2.18}$$

Let $a(n,k) = [x^{n-k}](1+(1+t)x+tx^2)^n$. Since $\left(1 + \frac{1+t}{\sqrt{t}}x + x^2\right)^n$ is palindromic we have

$[x^{2n-j}](1+(1+t)x+tx^2)^n = t^{n-j}[x^j](1+(1+t)x+tx^2)^n$ and thus

$[x^n](1+(1+t)x+tx^2)^{n-1} = t[x^{n-2}](1+(1+t)x+tx^2)^{n-1}$.

For $k \geq 1$ we have

$a(n,k) = [x^{n-k}](1+(1+t)x+tx^2)^n = [x^{n-k}](1+(1+t)x+tx^2)(1+(1+t)x+tx^2)^{n-1}$

$= [x^{n-1-(k-1)}](1+(1+t)x+tx^2)^{n-1} + (1+t)[x^{n-1-k}](1+(1+t)x+tx^2)^{n-1} + t[x^{n-1-(k+1)}](1+(1+t)x+tx^2)^{n-1}$

$= a(n-1, k-1) + (1+t)a(n-1, k) + ta(n-1, k+1)$.

For $k = 0$ we get

$a(n,0) = [x^n](1+(1+t)x+tx^2)^n$

$= [x^n](1+(1+t)x+tx^2)^{n-1} + (1+t)[x^{n-1-0}](1+(1+t)x+tx^2)^{n-1} + 2t[x^{n-1-(1)}](1+(1+t)x+tx^2)^{n-1}$

$= ta(n-1, k-1) + (1+t)a(n-1, 0) + ta(n-1, 1) = (1+t)a(n-1, 0) + 2ta(n-1, 1)$.

Another formula for $n > 0$ is

$$D_{n,k}(t) = \sum_{j=0}^{n} \binom{n}{j}\binom{n}{k+j} t^j. \qquad (2.19)$$

This follows from $(1+x+tx(1+x))^n = \sum_{j=0}^{n}\binom{n}{j} t^j x^j (1+x)^{n-j}$ by considering the coefficient of $x^{n-k}$.

By (2.17) the polynomials $D_{n,k}(t)$ can also been interpreted as the weight of all NSEW-paths of length $n$ and whose endpoint is on height $k$ with weights $w(E) = w(N) = 1$, $w(W) = t$, $w(S) = 2t$ if the endpoint of $S$ is on the $x$-axis and $w(S) = t$ else.

Let for example $n = 2$ and $k = 0$. Then we have $w(EE) = 1$, $w(WW) = t^2$, $w(NS) = 2t$, $w(EW) = w(WE) = t$. For $n = 2$ and $k = 1$ we get $w(NE) = w(EN) = 1$ and $w(WN) = w(NW) = t$.

The first terms of the sequence $(D_{n,0}(t), D_{n,1}(t), \cdots, D_{n,n}(t))_{n \geq 0}$ are

```
1
1 + t                    1
1 + 4 t + t²             2 + 2 t            1
1 + 9 t + 9 t² + t³      3 + 9 t + 3 t²     3 + 3 t        1
1 + 16 t + 36 t² + 16 t³ + t⁴   4 + 24 t + 24 t² + 4 t³   6 + 16 t + 6 t²   4 + 4 t   1
```

For $t=1$ $D_{n,k}(t)$ reduces to $D_{n,k}(1) = \binom{2n}{n-k}$ and we get the triangle OEIS [12], A094527,

$$\begin{pmatrix} 1 & 0 & 0 & 0 & 0 \\ 2 & 1 & 0 & 0 & 0 \\ 6 & 4 & 1 & 0 & 0 \\ 20 & 15 & 6 & 1 & 0 \\ 70 & 56 & 28 & 8 & 1 \end{pmatrix}$$

For $t=2$ we get OEIS [12], A118384,

$$\begin{pmatrix} 1 & 0 & 0 & 0 & 0 \\ 3 & 1 & 0 & 0 & 0 \\ 13 & 6 & 1 & 0 & 0 \\ 63 & 33 & 9 & 1 & 0 \\ 321 & 180 & 62 & 12 & 1 \end{pmatrix}$$

The polynomials $D_{n,k}(t)$ are gamma-nonnegative. More precisely we have

$$D_{n,k}(t) = \sum_{j=0}^{\lfloor \frac{n-k}{2} \rfloor} \binom{2j+k}{j} \binom{n}{2j+k} t^j (1+t)^{n-k-2j}. \tag{2.20}$$

The proof is analogous to the corresponding proof of (1.26).

For each non-negative NSEW- path $u_1 \cdots u_n$ with $u_i \in \{N, S, E, W\}$ whose endpoint is on height $k$ there are $i$ terms $f(u_j)$ negative and $i+k$ terms $f(u_j) = 1$ for some $i$. We can choose $2i+k$ places where $u_j = N$ or $u_j = S$ in $\binom{n}{2i+k}$ ways. By (2.2) for $t=1$ the weight of all non-negative paths is $\binom{k+2i}{i}$. The remaining $n-2i-k$ places can arbitrarily be filled with $W$ or $E$. Therefore for arbitrary $t$ the weight of all such paths is
$\binom{n}{2i+k}\binom{k+2i}{i} t^i (1+t)^{n-k-2i}.$

Let $M_0$ be the linear functional defined by $M_0(R_n(x,t)) = [n=0]$. Then (2.16) and (2.19) imply

$$M_0(x^n) = M_n(t). \tag{2.21}$$

This result can be found in [1] and [13] and is implicitly contained in [17].

Formula (2.16) implies $x^{2n} = \sum_{k=0}^{n} D_{n,k}(t) L_{2k}(x,t)$ and therefore

$$M(x^{2n}) = D_{n,0}(t) = M_n(t). \tag{2.22}$$

In the same way there are rational functions $E_{n,k}(t)$ such that $x^{2n+1} = \sum_{k=0}^{n} E_{n,k}(t) L_{2k+1}(x,t)$ which implies $M(x^{2n+1}) = 0$. This gives

**Theorem 2 ([1], [13], [17])**

*Let $M$ be the linear functional defined by $M(L_n(x,t)) = [n=0]$. Then the moments satisfy*

$$\begin{aligned} M(x^{2n}) &= M_n(t), \\ M(x^{2n+1}) &= 0. \end{aligned} \tag{2.23}$$

Let us now compute the generating functions $f_k(z,t) = \sum_{n \geq 0} D_{n,k}(t) z^n$.

We get $f_k(z,t) = z\big(f_{k-1}(z,t) + (1+t)f_k(z,t) + t f_{k+1}(z,t)\big)$ for $k > 0$ and
$f_0(z,t) = 1 + (1+t)z f_0(z,t) + 2tz f_1(z,t)$.

This gives $f_k(z,t) = z^k f_0(z,t) f(z,t)^k$ with

$$f(z,t) = \frac{1-(1+t)z - \sqrt{(1-(1+t)z)^2 - 4tz^2}}{2tz^2} = \frac{C(t,z) - 1}{z} \quad \text{by (1.29). Thus}$$

$$f_0(z,t) = \frac{1}{1-(1+t)z - 2tz^2 f(z,t)} = \frac{1}{\sqrt{(1-(1+t)z)^2 - 4tz^2}}.$$

This gives

$$M(t,z) = \sum_{n \geq 0} M_n(t) z^n = \frac{1}{\sqrt{(1-(1+t)z)^2 - 4tz^2}} \tag{2.24}$$

and

$$\sum_{n \geq 0} D_{n,k}(t) z^n = M(t,z)\big(C(t,z) - 1\big)^k. \tag{2.25}$$

**Corollary**

*Let*

$$c_n(m,t) = \sum_{k=0}^{n-1}\binom{n-1}{k}\binom{n+m}{k+m}\frac{m}{n+m}t^k$$

*with $c_0(m,t)=1$ be the $m$-fold convolution of $C_n(t)$ with itself (cf. (3.2)).*

*Then for $m \geq 1$*

$$\frac{1}{\prod_{j=0}^{m-1}(n-j)}\sum_{k=0}^{n}\left(\frac{\partial^m}{\partial t^m}D_{n,k}(t)\right)R_k(x,t) = \sum_{j=0}^{n-m}c_{n-m-j}(m,t)x^j. \tag{2.26}$$

**Proof**

By (3.4) we have

$$\frac{\partial^m}{\partial t^m}\sum_{n\geq 0}\frac{D_{n+m,k}(t)}{(n+m)\cdots(n+1)}z^n = C(t,z)^m\sum_{n\geq 0}D_{n,k}(t)z^n.$$

Therefore the left-hand side of (2.26) is the coefficient of $z^{n-m}$ of the power series

$$C(t,z)^m\sum_{n\geq 0}\sum_{k=0}^{n}D_{n,k}(t)R_k(x,t)z^n = C(t,z)^m\sum_{n\geq 0}x^n z^n = \sum_{i\geq 0}c_i(m,t)z^i\sum_{\ell\geq 0}x^\ell z^\ell$$

and the coefficient of $z^{n-m}$ is $\sum_{j=0}^{n-m}c_{n-m-j}(m,t)x^j$.

Since $\frac{\partial^m}{\partial t^m}D_{n,k}(t)\Big|_{t=1} = \sum_{k=0}^{n}\binom{n}{j}\binom{n}{k+j}\binom{j}{m} = \binom{n}{m}\sum_{k=0}^{n}\binom{n}{k+j}\binom{n-m}{n-j} = \binom{n}{m}\binom{2n-m}{k+n}$

(2.26) for $t=1$ implies

$$\sum_{k=0}^{n-m}\binom{2n-m}{n+k}L_{2k}(x) = \sum_{j=0}^{n-m}c_j(m,1)x^{2(n-m-j)} = \sum_{j=0}^{n-m}\frac{m}{m+2j}\binom{m+2j}{j}x^{2(n-m-j)}.$$

For $m=1$ this reduces to

$$\sum_{k=0}^{n-1}\binom{2n-1}{n+k}L_{2k}(x) = \sum_{j=0}^{n-1}\frac{1}{1+2j}\binom{1+2j}{j}x^{2(n-1-j)} = \sum_{j=0}^{n-1}C_j x^{2(n-1-j)}.$$

It seems that there are also similar extensions of (1.22) and (1.33).

**Conjecture 1**

$$\sum_{k=0}^{n}\left(\frac{\partial^m}{\partial t^m} A_{n,k}(t)\right) P_k(x,t) = \prod_{j=1}^{m-1}(n-j) \sum_{j=0}^{n-m-1}(j+1)x^{j+1}c_{n-m-j-1}(m,t), \qquad (2.27)$$

$$\sum_{k=0}^{n}\left(\frac{\partial^m}{\partial t^m} B_{n,k}(t)\right) Q_k(x,t) = \prod_{j=1}^{m-1}(n+1-j) \sum_{j=0}^{n-m}(j+1)x^{j}c_{n-m-j}(m,t). \qquad (2.28)$$

Let me only mention one special case for $m=1$.

Since $\left.\frac{\partial B_n(k,t)}{\partial t}\right|_{t=1} = (k+1)\binom{2n+1}{n-k-1}$ we get

$$\sum_{k=0}^{n}(k+1)\binom{2n+1}{n-k-1} F_{2k+1}(x) = \sum_{j=0}^{n-1}(j+1)C_{n-1-j}x^{2j+1}.$$

**2.3. The polynomials $S_n(x,t) = \dfrac{L_{2n+1}(\sqrt{x},t)}{\sqrt{x}}$.**

Let $\sigma_0(t) = \dfrac{1+4t+t^2}{1+t}$ and $\sigma_n(t) = \dfrac{1+t^{n+1}}{1+t^n} + t\dfrac{1+t^n}{1+t^{n+1}}$.

The polynomials

$$S_n(x,t) = \frac{L_{2n+1}(\sqrt{x},t)}{\sqrt{x}} \qquad (2.29)$$

satisfy the recursion

$$S_n(x,t) = (x - \sigma(n-1,t))S_{n-1}(x,t) - \frac{t(1+t^{n-2})(1+t^n)}{(1+t^{n-1})^2} S_{n-2}(x,t)$$

with initial values $S_0(x,t) = 1$ and $S_1(x,t) = x - \dfrac{1+4t+t^2}{1+t}$.

**Theorem 3**

*The polynomials $S_n(x,t)$ are explicitly given by*

$$S_n(x,t) = \frac{1}{1+t^n}\sum_{k=0}^{n}(-1)^{n-k}G_{n,k}(t)x^k \qquad (2.30)$$

*with*

$$G_{n,k}(t) = \sum_{j=0}^{n-k} \frac{\binom{j+k}{k}\binom{n-j-1}{k-1}(n(k+1)-j)}{k(k+1)}\left(t^j + t^{2n-k-j}\right). \quad (2.31)$$

*for $k > 0$ and*

$$G_{n,0}(t) = (2n+1)t^n + \sum_{j=0}^{2n} t^j. \quad (2.32)$$

The first terms of the sequence $\left(G_{n,0}(t), G_{n,1}(t), \cdots, G_{n,n}(t)\right)_{n \geq 0}$ are

```
2
1 + 4 t + t²
1 + t + 6 t² + t³ + t⁴                 1 + t
1 + t + t² + 8 t³ + t⁴ + t⁵ + t⁶       2 + 3 t + 3 t² + 2 t³                              1 + t²
                                        3 + 5 t + 6 t² + 6 t³ + 5 t⁴ + 3 t⁵    3 + 4 t + 4 t³ + 3 t⁴    1 + t³
```

To prove this observe that by (2.10) we get

$$xS_n(x,t) = R_{n+1}(x,t) + \tau(2n,t)R_n(x,t).$$

This is equivalent with

$$\left[x^{k+1}\right]\left((1+t^n)R_{n+1}(x,t) + (1+t^{n+1})R_n(x,t)\right) = (-1)^{n-k} G_{n,k}(t).$$

Let us first consider the coefficient of $t^j$ with $j < n$.

Comparing coefficients gives the easily verified identity

$$-\binom{n}{k+1}\binom{n-k-1}{j}\frac{\binom{k+j}{j}}{\binom{n-1}{j}} + \binom{n+1}{k+1}\binom{n-k}{j}\frac{\binom{k+j}{j}}{\binom{n}{j}} =$$

$$\frac{\binom{n-j-1}{k-1}\binom{k+j}{j}((k+1)n-j)}{k(k+1)}.$$

Now let us consider the coefficient of $t^{2n-k-j}$. Here we have to show that

$$(-1)^{n-k}[t^{n-k-j}x^{k+1}]\left(R_{n+1}(x,t) + tR_n(x,t)\right) = \frac{\binom{j+k}{k}\binom{n-j-1}{k-1}(n(k+1)-j)}{k(k+1)}.$$

The left-hand side is

$$\binom{n+1}{k+1}\binom{n-k}{j}\binom{n-j}{k}\frac{1}{\binom{n}{k+j}} - \binom{n}{k+1}\binom{n-k-1}{j-1}\binom{n-j}{k}\frac{1}{\binom{n-1}{k+j-1}}$$

which can be simplified to give the right-hand side.

The coefficients of $G_{n,k}(t)$ are related to the numbers $g(n,j,k)$ in OEIS [12] A051340, A141419, A185874, A185875, A185876.

**Theorem 4**

*The functions $E_{n,k}(t)$ which satisfy*

$$\sum_{k=0}^{n} E_{n,k}(t) S_k(x,t) = x^n \tag{2.33}$$

*are*

$$E_{n,k}(t) = \frac{\sum_{j=0}^{n-k} \binom{n}{k+j}\binom{n+1}{j}(t^j + t^{n+1-j})}{1+t^{k+1}} \tag{2.34}$$

*for $n \geq k$ and $E_{n,k}(t) = 0$ else.*

*As special case note that*

$$E_{n,0}(t) = \frac{\sum_{j=0}^{n} \binom{n}{j}\binom{n+1}{j}(t^j + t^{n+1-j})}{1+t} = \frac{\sum_{j=0}^{n+1} \binom{n+1}{j}^2 t^j}{1+t} = \frac{M_{n+1}(t)}{1+t}. \tag{2.35}$$

**Proof**

By (1.1) this follows from

$$E_{n,k}(t) = D_{n,k}(t) + \tau(2k+1) D_{n,k+1}(t) = \sum_{j=0}^{n-k} \binom{n}{j}\binom{n}{k+j} t^j + \frac{t(1+t^k)}{1+t^{k+1}} \sum_{j=0}^{n} \binom{n}{j}\binom{n}{k+j+1} t^j$$

$$= \frac{1}{1+t^{k+1}}\left( \sum_{j=0}^{n-k}\binom{n+1}{j}\binom{n}{k+j} t^j + \sum_{j=0}^{n-k}\binom{n}{j}\binom{n+1}{k+j+1} t^{j+k+1} \right)$$

$$= \frac{1}{1+t^{k+1}}\left( \sum_{j=0}^{n-k}\binom{n}{j}\binom{n}{k+j} t^j + \sum_{j=0}^{n-k}\binom{n}{k+j}\binom{n+1}{j} t^{n-j+1} \right).$$

Thus the linear functional $M_1$ defined by $M_1(S_n(x,t)) = [n = 0]$ has the moments

$$M_1(x^n) = \frac{M_{n+1}(t)}{1+t}. \tag{2.36}$$

The first terms of the triangle $\left((1+t)E_{n,0}(t), (1+t^2)E_{n,1}(t), \cdots, (1+t^{n+1})E_{n,n}(t)\right)_{n \geq 0}$ are

```
1 + t
1 + 4 t + t²
1 + 9 t + 9 t² + t³                1 + t²
1 + 16 t + 36 t² + 16 t³ + t⁴      2 + 3 t + 3 t² + 2 t³              1 + t³
                                   3 + 12 t + 12 t² + 12 t³ + 3 t⁴   3 + 4 t + 4 t³ + 3 t⁴    1 + t⁴
```

The first terms of the triangle $\left(E_{n,0}(2), E_{n,1}(2), \cdots, E_{n,n}(2)\right)_{n \geq 0}$ are

$$\begin{pmatrix}
1 & 0 & 0 & 0 & 0 & 0 & 0 & 0 \\
\frac{13}{3} & 1 & 0 & 0 & 0 & 0 & 0 & 0 \\
21 & \frac{36}{5} & 1 & 0 & 0 & 0 & 0 & 0 \\
107 & \frac{219}{5} & \frac{91}{9} & 1 & 0 & 0 & 0 & 0 \\
561 & \frac{1272}{5} & \frac{226}{3} & \frac{222}{17} & 1 & 0 & 0 & 0 \\
\frac{8989}{3} & 1453 & \frac{4510}{9} & \frac{1970}{17} & \frac{529}{33} & 1 & 0 & 0 \\
16213 & 8244 & 3155 & \frac{14886}{17} & \frac{1821}{11} & \frac{1236}{65} & 1 & 0 \\
\frac{265729}{3} & \frac{233303}{5} & \frac{57799}{3} & \frac{103299}{17} & \frac{46403}{33} & \frac{14581}{65} & \frac{2839}{129} & 1
\end{pmatrix}$$

Note that the first column contains the numbers $E_{n,0}(2) = \frac{M_{n+1}(2)}{3}$. By [7], Theorem 5.8, the Delannoy numbers $M_n(2)$ are multiples of 3, i.e. $E_{n-1,0}(2) \in \mathbb{N}$, if and only if the base 3 representation of $n$ contains at least one 1. This is sequence OEIS [12], A081606, $(1, 3, 4, 5, 7, 9, \cdots)$.

### 3. Convolutions of Narayana polynomials.

Finally we want to derive some convolution formulae. By (1.41) we have

$$C(t, z) = \sum_{n \geq 0} C_n(t) z^n = \frac{1 + z(t-1) - \sqrt{1 - 2z(t+1) + z^2(t-1)^2}}{2tz}$$

or equivalently

$$tzC(t,z)^2 = C(t,z) - 1 - zC(t,z) + tzC(t,z). \tag{3.1}$$

We will show that

$$C(t,z)^m = \sum_{n \geq 0} c_n(m,t) z^n \tag{3.2}$$

with

$$c_n(m,t) = \sum_{k=0}^{n-1} \binom{n-1}{k}\binom{n+m}{k+m}\frac{m}{n+m} t^k \qquad (3.3)$$

and $c_0(m,t) = 1$.

Note that $c_n(1,t) = \sum_{k=0}^{n-1}\binom{n-1}{k}\binom{n+1}{k+1}\frac{1}{n+1}t^k = \sum_{k=0}^{n-1}\binom{n-1}{k}\binom{n}{k}\frac{1}{k+1}t^k = C_n(t)$.

It suffices to show that

$$tzC(t,z)^m = C(t,z)^{m-1}(1+z(t-1)) - C(t,z)^{m-2}$$

holds if we replace $C(t,z)^m$ by $\sum_{n\geq 0} c_n(m,t)z^n$.

The coefficient of $z^{n+1}$ is

$$tc_n(m,t) = c_{n+1}(m-1,t) + (t-1)c_n(m-1,t) - c_{n+1}(m-2,t).$$

The coefficient of $t^{k+1}$ is

$$\binom{n-1}{k}\binom{n+m}{k+m}\frac{m}{n+m} = \binom{n}{k+1}\binom{n+m}{k+m}\frac{m-1}{n+m} + \binom{n-1}{k}\binom{n+m-1}{k+m-1}\frac{m-1}{n+m-1}$$

$$-\binom{n-1}{k+1}\binom{n+m-1}{k+m}\frac{m-1}{n+m-1} - \binom{n}{k+1}\binom{n+m-1}{k+m-1}\frac{m-2}{n+m-1}$$

Dividing by $\binom{n-1}{k}\binom{n+m-1}{k+m-1}$ this gives

$$\frac{m}{k+m} = \frac{n}{k+1}\frac{m-1}{k+m} + \frac{m-1}{n+m-1} - \frac{n-k-1}{k+1}\frac{n-k}{k+m}\frac{m-1}{n+m-1} - \frac{n}{k+1}\frac{m-2}{n+m-1}$$

which is easily verified.

More generally we want to show that

$$\frac{\partial^m}{\partial t^m}\sum_{n\geq 0}\frac{D_{n+m,k}(t)}{(n+m)\cdots(n+1)}z^n = C(t,z)^m \sum_{n\geq 0} D_{n,k}(t)z^n. \qquad (3.4)$$

The coefficient of $z^n$ of the left-hand side is

$$v(n,m,k) = \sum_{j=0}^{n}\frac{\binom{n+m}{j}\binom{n+m}{j+k}\binom{j}{m}}{\binom{n+m}{m}}t^{j-m}$$

As above it suffices to verify that

$$tzC(t,z)^m \sum_{n\geq 0} D_{n,k}(t)z^n = C(t,z)^{m-1} \sum_{n\geq 0} D_{n,k}(t)z^n \left(1+z(t-1)\right) - C(t,z)^{m-2} \sum_{n\geq 0} D_{n,k}(t)z^n$$

or

$$tv(n,m,k) = v(n+1,m-1,k) + (t-1)v(n,m-1,k) - v(n+1,m-2,k).$$

This can easily be verified.

For $t=1$ formula (3.2) reduces to the well-known formula

$$C(1,z)^m = \left(\frac{1-\sqrt{1-4z}}{2z}\right)^m = \sum_{n\geq 0} \frac{m}{2n+m} \binom{2n+m}{n} z^n. \tag{3.5}$$

A well-known convolution formula for the central binomial coefficients is

$$\sum_{k=0}^n \binom{2k}{k}\binom{2(n-k)}{n-k} = 4^n. \tag{3.6}$$

A computational proof follows immediately by squaring the generating function (2.4).

For the $m-$fold convolution we get

$$u_m(n) = \sum_{i_1+\cdots+i_m=n} \binom{2i_1}{i_1}\binom{2i_2}{i_2}\cdots\binom{2i_m}{i_m} = 4^n \binom{\frac{m}{2}+n-1}{n} \tag{3.7}$$

since

$$\left(\sum_{n\geq 0}\binom{2n}{n}x^n\right)^m = (1-4x)^{-\frac{m}{2}} = \sum_k \binom{-\frac{m}{2}}{k}(-4)^k x^k = \sum_k \binom{\frac{m}{2}+k-1}{k} 4^k x^k.$$

A combinatorial proof has been given in [8].

I want now to compute the corresponding convolutions of the polynomials $M_n(t)$.

Their generating function is

$$\sum_{n\geq 0} M_n(t)x^n = \frac{1}{\sqrt{(1+(1-t)x)^2 - 4x}}. \tag{3.8}$$

Let

$$\left(\frac{1}{\sqrt{(1+(1-t)x)^2 - 4x}}\right)^m = \sum_{n\geq 0} u_m(n,t)x^n. \tag{3.9}$$

Then we get

**Theorem 5**

$$u_m(n,t) = \sum_{k \geq 0} \binom{n+m-1}{m-1}\binom{n}{k} \frac{\prod_{j=0}^{k-1}(2n+m-1-2j)}{\prod_{j=0}^{k-1}(2k+m-1-2j)} t^k. \quad (3.10)$$

To prove these identities by induction observe that

$$u_{m-2}(n,t) = u_m(n,t) - (1+t)u_m(n-1,t) + (1-t)^2 u_m(n-2,t)$$

holds for all $n$.

The first 5 terms of $u_1(n,t), u_2(n,t), \cdots, u_5(n,t)$ are

```
1   1 + t    1 + 4 t + t²     1 + 9 t + 9 t² + t³       1 + 16 t + 36 t² + 16 t³ + t⁴
1   2 + 2 t  3 + 10 t + 3 t²  4 + 28 t + 28 t² + 4 t³   5 + 60 t + 126 t² + 60 t³ + 5 t⁴
1   3 + 3 t  6 + 18 t + 6 t²  10 + 60 t + 60 t² + 10 t³ 15 + 150 t + 300 t² + 150 t³ + 15 t⁴
1   4 + 4 t  10 + 28 t + 10 t² 20 + 108 t + 108 t² + 20 t³ 35 + 308 t + 594 t² + 308 t³ + 35 t⁴
1   5 + 5 t  15 + 40 t + 15 t² 35 + 175 t + 175 t² + 35 t³ 70 + 560 t + 1050 t² + 560 t³ + 70 t⁴
```

All these polynomials are palindromic and gamma-nonnegative:

$$u_m(n,t) = \sum_{k=0}^{n} \binom{n+m-1}{m-1}\binom{2k}{k}\binom{n}{2k} \frac{(2k)!!}{\prod_{i=0}^{k-1}(m+2i+1)} t^k (1+t)^{n-2k}. \quad (3.11)$$

For the proof we make use of Gauss's theorem for hypergeometric polynomials

$$_2F_1\left(\begin{matrix}a,b\\c\end{matrix};1\right) = \frac{\Gamma(c)\Gamma(c-a-b)}{\Gamma(c-a)\Gamma(c-b)}. \quad (3.12)$$

By comparing coefficients of $t^k$ in (3.10) and (3.11) it suffices to show that

$$\sum_{j=0}^{k} \frac{\binom{2j}{j}\binom{n}{2j}(2j)!!\binom{n-2j}{k-j}}{\binom{n}{k}\prod_{i=0}^{j-1}(m+2i+1)} = \frac{\prod_{j=0}^{k-1}(2n+m-1-2j)}{\prod_{j=0}^{k-1}(2k+m-1-2j)}.$$

The left-hand side can we written as $\,_2F_1\!\left(\begin{matrix}-k,k-n\\ \frac{m+1}{2}\end{matrix};1\right)$ which by Gauss's Theorem equals

$$\frac{\Gamma\!\left(\frac{m+1}{2}\right)\Gamma\!\left(\frac{m+1}{2}+n\right)}{\Gamma\!\left(\frac{m+1}{2}+k\right)\Gamma\!\left(\frac{m+1}{2}+n-k\right)} = \frac{\prod_{j=0}^{k-1}(2n+m-1-2j)}{\prod_{j=0}^{k-1}(2k+m-1-2j)}.$$

Let us finally consider two special cases in detail.

For $m=2$ we get

$$u_2(n,t) = \sum_{k=0}^{n} M_k(t) M_{n-k}(t) = \frac{1}{2}\sum_{k=0}^{n}\binom{2n+2}{2k+1}t^k = \sum_{k}\binom{n+1}{2k}t^k \sum_{k}\binom{n+1}{2k+1}t^k. \quad (3.13)$$

For the generating function of $u_2(n,t^2)$ is

$$\sum_{n\geq 0} u_2(n,t^2) x^n = \frac{1}{(1+(1-t^2)x)^2 - 4x} = \frac{1}{4t}\left(\frac{(1+t)^2}{1-(1+t)^2 x} - \frac{(1-t)^2}{1-(1-t)^2 x}\right).$$

This implies

$$u_2(n,t^2) = \frac{(1+t)^{2n+2} - (1-t)^{2n+2}}{4t} = \frac{1}{2}\sum_{k=0}^{n}\binom{2n+2}{2k+1}t^{2k}.$$

The right-hand side follows from $(1+t)^{2n} - (1-t)^{2n} = \big((1+t)^n + (1-t)^n\big)\big((1+t)^n - (1-t)^n\big)$.

For $m=3$ we get

$$u_3(n,t) = \sum_{k}\binom{n+2}{2}\binom{n}{k}\frac{\binom{n+1}{k}}{\binom{k+1}{1}} t^k = \binom{n+2}{2}\sum_{k}\binom{n}{k}\binom{n+1}{k}\frac{1}{k+1}t^k = \binom{n+2}{2}C_{n+1}(t).$$

It would be interesting to find combinatorial interpretations of these results.